# Tables in Leibniz
# a challenge for the digital humanities

**Ariles Remaki**
ERC Philiumm (101020985), Sphere, CNRS, Université Paris Cité

## Introduction

Leibniz's image as a *master-builder of mathematical notations*[1] tends to overshadow the commentary on those aspects of his practice that are outside formalism. Yet, in a philosophical text entitled by the editors "*De Arte characteristica inventoriaque analytica combinatoriave in mathesi univsersali*" (A VI, 4, 315), probably written between 1679 and 1680, Leibniz asserts that characters and symbols alone do not constitute the essential tools of the art of invention. Alongside them stands another, much less commented tool: tables.

> Just as we need an Instrument to draw a circle precisely, by which the hand is guided, and this all the more so as we are less practiced in it, so to think straight we lack some sensitive Instruments, which I reduce to two main chapters: *Characters* and *Tables*. The former we lack for Analysis, the latter for Combination.[2]

But if Leibniz's reflections on symbols can be inscribed in the long historical development of his project of a universal characteristic, what can be said of this attitude towards tables? How did Leibniz come to think of tabular practice as one of the two main pillars of the art of invention?
In considering this question, we are led to focus our attention on the German philosopher's youth, and more particularly on the Parisian stay during which he learned the mathematics of his time. The numerous tables revealed in these documents not only shed light on the way in which Leibnizian conceptions of the *ars inveniendi* were shaped, but they also provide particularly interesting cases for questioning the notion of the table as a writing tool, a support for reasoning or even a scriptural practice.
However, the present contribution does not aim to address the origin of Leibnizian tabular practice head-on[3], but rather to show how the semantic, philological and anthropological problems raised by the study of the tables of the Parisian period respond to the serious difficulties presented by the enterprise of editing these same tables, and particularly the digital edition. Indeed, a great diversity of data entry methods goes unnoticed in a paper edition, literally crushed by printing. On the other hand, the implementation and structure of data is a fundamental aspect of digital publishing that cannot be separated from the display process. From then on, these technical choices must be the product of profound historiographical reflection, and the tables produced by the young Leibniz constitute a real challenge for this purpose.
By the way, one of the objectives of the ERC Philiumm project is precisely to reflect on the digital edition of these complex elements from the Leibnizian corpus.

---

[1] Florian Cajori: "Leibniz, the Master-Builder of Mathematical Notations" in: *Isis* 7/3 (1925), pp. 412–429.

[2] Gottfried Wilhelm Leibniz: *Mathesis universalis: écrits sur la mathématique universelle*, ed. by David Rabouin, Vrin, 2018, p.82

[3] See my thesis: Ariles Remaki: *L'art combinatoire en tant qu'art d'inventer chez Leibniz, sur la période 1672-1680*, thesis defended at the Université Paris Cité, 2021

## I. Tabular practice in Leibniz: a Parisian story

In order to better understand Leibniz's conception of tables, we introduce two distinct approaches to tabular practice: the computational approach and the encyclopedic approach. The first consists in seeing the table as a computational tool. The different terms of the table are produced by an algorithmic procedure based on the structure of the table. The second is to see the table as an inventory that allows for an exhaustive listing of a set of terms. It seems possible to bring these two approaches together by noting that the inventory does involve a classification and indexing procedure which is ultimately similar in nature to the algorithmic procedure linked to the first approach. But the main difference between the two approaches lies in the epistemological priority given to terms. The terms of a table seen in a computational way are synthetic products of the tabular practice, whereas in the encyclopedic approach, it is the table itself that is produced by the analysis of the set of terms and the resulting ranking procedure. The main, highly theoretical statements on the role of tables within Leibniz's mature philosophical texts are rather encyclopedic in approach. The tables are to serve the cumulative character of knowledge and provide a clear presentation medium for forging the architectonics of general science. In his early work, however, it is the tables of calculation that play a central role. On this basis alone, it would be possible to envisage that these two approaches are part of two independent intellectual paths within the Leibnizian corpus, one more theoretical and philosophical, the other more practical and mathematical. But these two approaches are in fact linked by a major aspect of Leibnizian thought, namely the art of combination. Thus, the analysis of the early tables, which are constructed according to the computational approach and which constitute the main object of the present contribution, provides a certain light to the later discourses on the encyclopedic role of the tables.

### a. *Double entry tables and lists of lists.*

Already in the *Dissertatio de Arte Combinatoria* we find traces of an important practice of tables to guide the young philosopher's thinking. During his stay in Paris, this practice developed even further around a tool for analysing numerical progressions that Leibniz called the triangle of differences. This scheme consists of listing the terms of a progression, then determining the progression of the differences of the neighbouring terms, then repeating the operation a certain number of times (in the example below, the triangle of differences applies to a geometric progression of ratio 3):

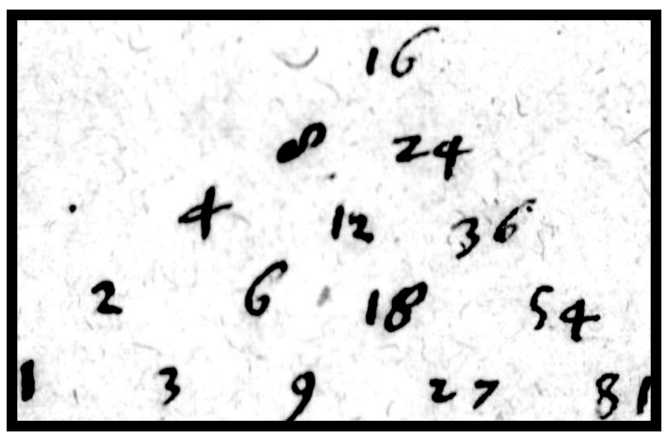

The triangle of differences is constructed as an ordered succession of progressions, which are themselves, by definition, ordered successions of terms. The internal structure of this scheme is therefore hierarchical: a term is first considered as an element of a progression, and this progression itself is part of the succession of sequences that form the triangle. We call this way of conceiving the table the "list of lists" strategy. Many elements in the *Dissertatio de Arte Combinatoria* testify to such an approach. On the one hand, some of the tables are entirely linear, i.e. the internal axis within the progressions is the same as that guiding the succession of progressions. On the other hand, for those tables that are truly two-dimensional, Leibniz explains their structure through the image of the generation of the rectangle by the translation of one orthogonal segment along another. This image again shows the hierarchical conception of the tabular structure: a fixed segment serves as a support for another segment, which is mobile.

The first works of the Parisian stay show a progressive evolution of the young Leibniz's practice on this point. Gradually, Leibniz sets up strategies that reverse the hierarchy within a table. In the end, the term becomes an element directly derived from the table, and the progressions are structures designed *a posteriori* from the terms. It is no longer the table that is synthetically derived from a succession of progressions, but rather the progressions that are analytically derived from a partition of the table. We choose to call this conception the "double entry table" strategy.

The distinction between these two approaches forms the basis of a very serious philological problem. Indeed, the two conceptions, "list of lists" and "double entry table" produce exactly the same table. It is therefore not possible to determine the position of the author by giving only the arrangement of terms in a table. Of course, there are many peritextual elements that could support the hypothesis that he chose one of the two approaches rather than the other. But this would suggest a psychologising approach which we do not subscribe to, and that is why we rather speak of strategy instead. This distinction does not therefore involve the psychological capacity to perceive the table according to this or that structure, because it is clear that, as early as the *Dissertatio de Arte Combinatoria*, Leibniz can consider horizontal, vertical or even diagonal progressions in the same table. What is essential, in our view, is the possibility of making this change of point of view operate, i.e. a concrete procedural practice that makes one point of view interact with the other and produce new results by confronting and articulating those distincts ways of reading. The "list of lists" strategy is characterised by the absence of these interactions, whereas the "double entry table" strategy fully integrates the interdependence of these different reading paths in the table.

### b. *Diachronic construction of a table*

In order for the study of Leibniz's drafts to serve the purpose of analysing his practice and the strategies he uses to discover new truths, it is natural to seek to understand the way in which these texts were written, and in particular the chronological order of the different elements that make up the manuscripts. This need is all the stronger in the case of the tables, for the reasons we have indicated about the different strategies for their design. However, the naively realistic attitude of wanting to re-establish the writing movement as it was performed by Leibniz seems to us to be as vain as the psychologising objective of re-establishing the reality of his cognitive processes. Here again, it is the material structure of the support that is the object of our analysis. Indeed, it should be remembered that the same manuscript element may have been written in an infinite number of different ways and that the same writing movement may result from an infinite number of distinct cognitive processes. Nevertheless, this does not mean that the historian must resign himself to considering the product of writing as a synchronic whole in which there is no trace of a dynamic practice distributed over time. We would simply point out that the hypotheses that the genetic analysis of texts makes it possible to support are always subject to the principle of parsimony and for this reason always include an irreducible element of retrospective reconstruction.
Let us take as an example a text by Leibniz regularly studied by historians of mathematics[4], namely *De numero jactuum in tesseris*[5], dated January 1676, in which Leibniz introduces the table of combinatorial numbers at the end of an exploration of a probability problem on dice throwing. On the back of the second folio are two versions of this table, written opposite each other:

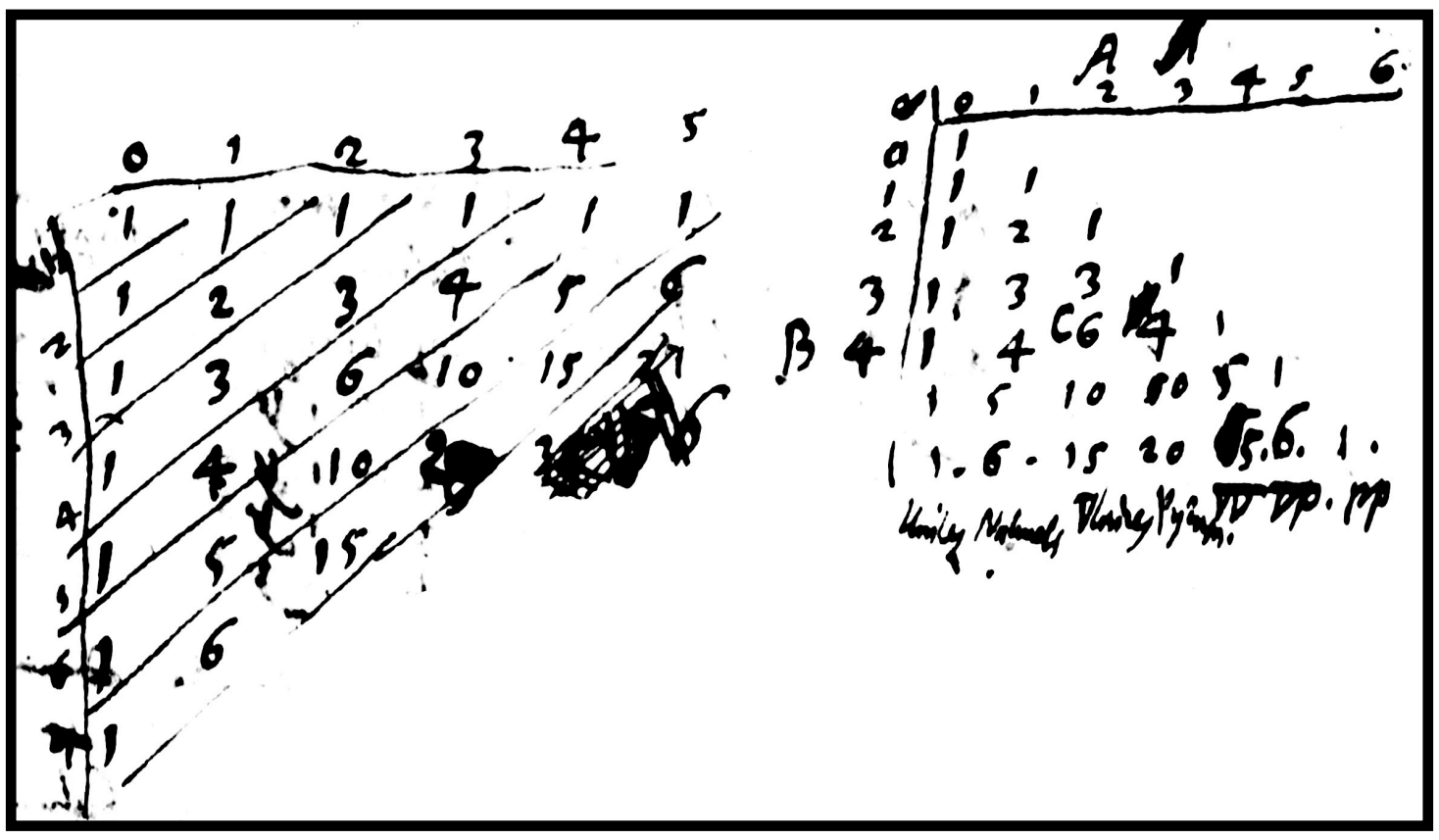

The table on the left corresponds to the one found in Pascal's *Traité du triangle arithmétique.* No explicit references are found in this text, but they are found attached to a very similar version of this table in another text, the *De Intercalatione numerorum combinatoriorum*[6] dated January 1675. The table on the right has a structure closer to the table of combinatorial numbers in the *Dissertatio de Arte Combinatoria*, published a year after the posthumous publication of Pascal's treatise, but without Leibniz having had access to it. The structure is similar, but there are substantial differences, notably the reversal of rows and columns and the removal of null terms.

---

[4] Gottfried Wilhelm Leibniz: *L'estime des apparences: 21 manuscrits de Leibniz sur les probabilités, la théorie des jeux, l'espérance de vie*, ed. by Marc Parmentier, Vrin, 1995, here pp.73-103; Eberhard Knobloch: *Die Mathematischen Studien von G. W. Leibniz zur Kombinatorik: auf Grund fast ausschliesslich handschriftlicher Aufzeichnungen dargelegt und kommentiert*, Studia liebnitiana supplementa 11. Wiesbaden: F. Steiner, 1973, here p.71; Dietrich Mahnke: "Leibniz auf der Suche nach einer allgemeinen Primzahlgleichung" in: *Bibliotheca mathematica* 13/3 (1912): pp. 29-61, here p.40.

[5] This manuscript is marked LH 35, 3 B 14, Bl.1-2. The part that interests us has been separated by the editors into another text, which is published in A VII 1, pp. 583-586 under the title that is not by Leibniz: *De numeris figuratis divisoribusque potestatum*.

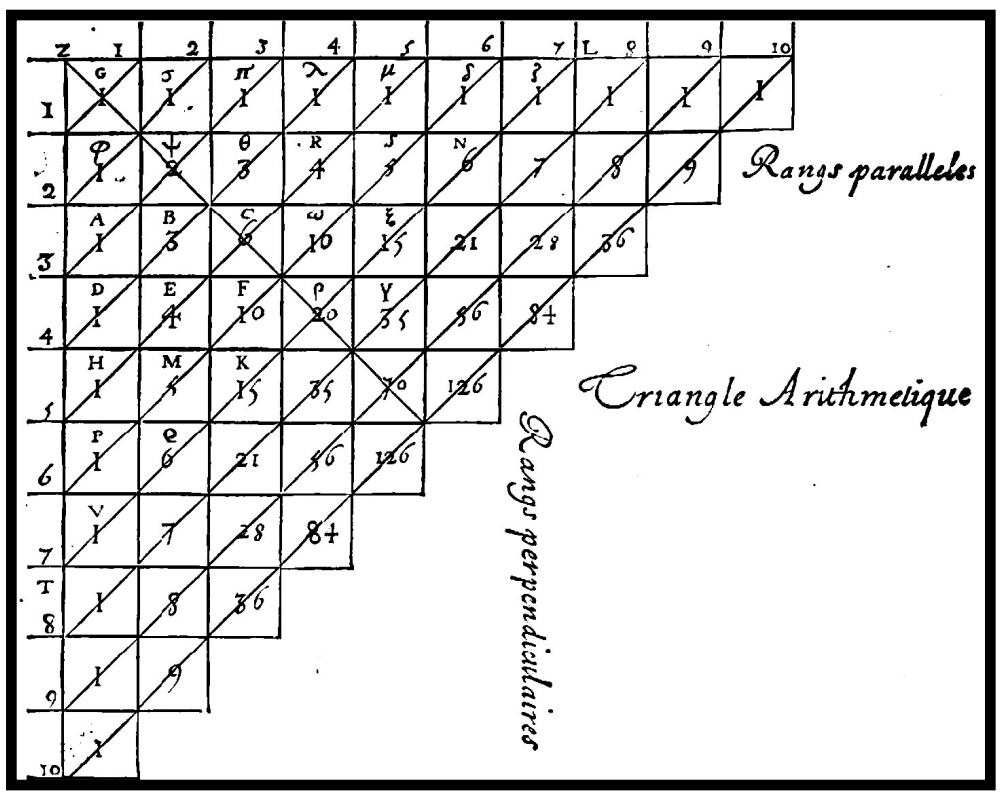

Pascal (1665)

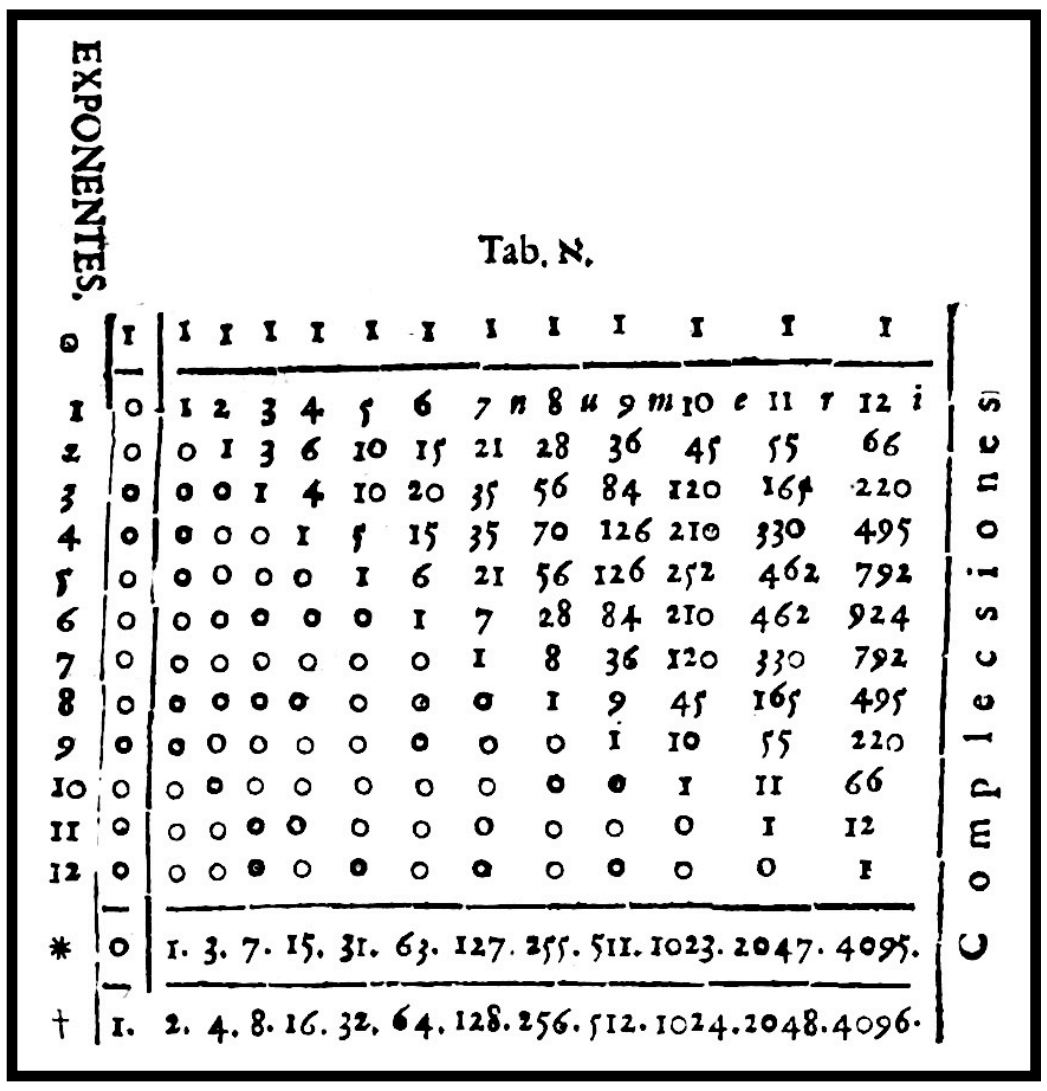

Leibniz (1666)

In the January 1676 text, Leibniz refers to the table on the right as the "arrangement [*dispositio*] of the arithmetic triangle for combinations". Indeed, the terms in Pascal's table are abstract coefficients that can be interpreted *a posteriori* as quantities that count combinations, but this is not how they are introduced or defined. In contrast, in the *Dissertatio de Arte Combinatoria*, the terms of the table are indeed defined as "complexions", i.e. a number of combinations. From this point on, in the text of 1676, Leibniz changes his object of study: he leaves the combinatorial problem of dice throwing to focus on an arithmetical problem of divisibility. In this second problem, the terms of the table are no longer interpreted as combination numbers, but as figured numbers. Now the figured numbers are also indexed in the table on the right, although they are arranged in a strange way, since each vertical progression does not start at the same height. The question we wish to address, through this example, is the following: are the two tables Leibniz writes face to face, two representations for the two interpretations (the left-hand table for the figured numbers and the right-hand table for the combinations)? Or is Pascal's table a copy that represents the general table, while the left-hand table would then represent both interpretations? As we have seen, this second hypothesis is based on the double indexing of the right-hand table, which thus provides elements for both the combinatorial and the geometric interpretation of the terms, while the indexing of the left-hand table is abstract. However, a genetic analysis of the table provides arguments for the first hypothesis. Indeed, in the left-hand table, under the erasures, one can clearly read the missing terms of the third and fourth line, which indicates that Leibniz constructed the left-hand table in the same way as the right-hand table, i.e. according to horizontal progressions. In the left-hand table, however, the horizontal progressions are those of the figured numbers. Moreover, still under the erasures, we can read the index "(y)" between the last two terms, "35" and "56" of the fourth line. And *y* is defined hereafter as the side of the figured numbers.

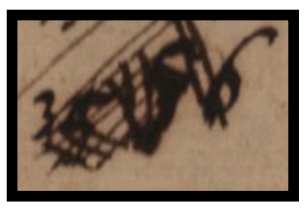

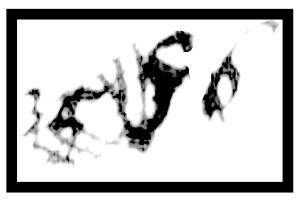

These palaeographic elements therefore support the hypothesis that the left-hand table, i.e. Pascal's, represents the figured numbers while the right-hand table represents the combinations. The consequences of this hypothesis are that a change of interpretation, a change of strategy implies a change of layout. Whereas the second hypothesis, which combines the two interpretations in the single right-hand table, implies that, on the contrary, the same layout can be associated with different strategies and interpretations.

---

[6] Gottfried Wilhelm Leibniz: *Die Mathematischen Studien von G. W. Leibniz zur Kombinatorik.* ed. by Eberhard Knobloch. Studia liebnitiana 16. Wiesbaden: F. Steiner, 1976, here pp.16-36.

These details might seem anecdotal. But they actually engage an interpretative dilemma that is explicitly embodied in the editorial choices. Thus Marc Parmentier[7], who seems to support the first hypothesis, chooses to edit the table by including the terms that are under the erasures (but not the index "(y)"). He removes the diagonal lines and thus obtains a table whose structure highlights the horizontal progressions, i.e. the figured numbers:

| | 0 | 1 | 2 | 3 | 4 | 5 |
|---|---|---|---|---|---|---|
| 1 | 1 | 1 | 1 | 1 | 1 | 1 |
| 2 | 1 | 2 | 3 | 4 | 5 | 6 |
| 3 | 1 | 3 | 6 | 10 | 15 | 21 |
| 4 | 1 | 4 | 10 | 20 | 35 | 56 |
| 5 | 1 | 5 | 15 | | | |
| 6 | 1 | 6 | | | | |
| 7 | 1 | | | | | |

| | 0 | 1 | A<br>2 | 3 | 4 | 5 | 6 |
|---|---|---|---|---|---|---|---|
| 0 | 1 | | | | | | |
| 1 | 1 | 1 | | | | | |
| 2 | 1 | 2 | 1 | | | | |
| 3 | 1 | 3 | 3 | 1 | | | |
| B4 | 1 | 4 | C6 | 4 | 1 | | |
| | 1 | 5 | 10 | 10 | 5 | 1 | |
| | 1 | 6 | 15 | 20 | 15 | 6 | 1 |
| | uni<br>tés | Natu<br>rels | Tri<br>angl | Pyra<br>mid | TT | Tp | pp[24] |

On the contrary, Eberhard Knobloch opts for an edition[8] that is more faithful to the manuscript, and completely removes the elements that have been erased[9] .

| | 0 | 1 | 2 | 3 | 4 | 5 |
|---|---|---|---|---|---|---|
| | 1 | 1 | 1 | 1 | 1 | 1 |
| 1 | 1 | 2 | 3 | 4 | 5 | 6 |
| 2 | 1 | 3 | 6 | 10 | 15 | |
| 3 | 1 | 4 | 10 | 20 | | |
| 4 | 1 | 5 | 15 | | | |
| 5 | 1 | 6 | | | | |
| 6 | 1 | | | | | |
| 7 | | | | | | |

| | 0 | 1 | A<br>2 | 3 | 4 | 5 | 6 |
|---|---|---|---|---|---|---|---|
| 0 | 1 | | | | | | |
| 1 | 1 | 1 | | | | | |
| 2 | 1 | 2 | 1 | | | | |
| 3 | 1 | 3 | 3 | 1 | | | |
| B 4 | 1 | 4 C | 6 | 4 | 1 | | |
| | 1 | 5 | 10 | 10 | 5 | 1 | |
| | 1 | 6 | 15 | 20 | 15 | 6 | 1 |
| | Unitéz | Naturels | ▽ laires | Pyram. | ▽▽ | ▽ p. | Pp. |

This choice is certainly more neutral than Parmentier's, but makes the second hypothesis much more credible than the first.

These two editions have chosen to avoid the genetic aspect, which is fundamental to the critical construction of the text. The erasures actually support a third hypothesis based on a genetic approach and which would suppose that Leibniz first opted for a strategy in line with the first hypothesis, i.e. two tables for two interpretations, and then he would have changed his point of view and placed both interpretations on the right-hand table.

This example does not only illustrate the contribution of the genetic approach to the exegesis of the young Leibniz's manuscript tables. It shows how the genetic approach must contribute to the undertaking of critical editing. This presupposes the integration of the sources into a diachronic narrative whose structure is based on material, palaeographic and philological elements. The contribution of a temporal dimension conflicts with the synchronic character of the paper edition, as we have already pointed out, which partly explains the editing choices made by Knobloch and Parmentier in the example we have given. This is why digital publishing is a great opportunity to address these issues.

---

[7] Gottfried Wilhelm Leibniz: *L'estime des apparences: 21 manuscrits de Leibniz sur les probabilités, la théorie des jeux, l'espérance de vie*, ed. by Marc Parmentier, Vrin, 1995 - p.80

[8] Eberhard Knobloch: *Die Mathematischen Studien von G. W. Leibniz zur Kombinatorik: auf Grund fast ausschliesslich handschriftlicher Aufzeichnungen dargelegt und kommentiert*, Studia liebnitiana supplementa 11.Wiesbaden: F. Steiner, 1973 - p.71

[9] This edition of the tables is reproduced, almost identically, in volume A VII, 1 of the *Akademie*, a volume edited by E. Knobloch among others.

However, the computer implementation of these aspects also has serious limitations. Indeed, we have shown that the genetic approach commits us to consider the sources as the product of material practices, whereas a digital edition forces us to express them through a program, which can always be reduced to a string of characters. This difficulty gives us the opportunity to question the very nature of the table, as a data structure and a practical instrument.

## II. New semantic approaches: diagrammatology

The term "diagrammatology" was introduced in the early 1980s by William Mitchell[10], an art historian, in order to criticise a certain idealist and Platonic approach in the historiography of the time. Mitchell defended a more Aristotelian position, in which the diagram constitutes the necessary interface between an intelligible form and its material medium of expression. Diagrammatology therefore consists in the study of the relations that are played out on this border. Thus, the diagram constitutes an object of study where metaphysics, phenomenology and semiology converge, as presented by the semiologist Frederik Stjernfelt[11], taking the term from Mitchell. Other, more anthropological approaches have been developed by historians of science, notably Kenneth Manders[12] through the notion of diagram control, a notion that is similar to that of the anthropologist Edwin Hutchins' material anchors[13]. These works, which bear witness to a profound change in the history of ideas about representations over the last forty years, agree on this point: the diagram necessarily brings the intelligible into dialogue with the sensible.

### *a. Is it all text or is it all image?*

The definition and characterisation of diagrams is a very topical issue. At first sight, the philosophical and philological approach to this question lies at the heart of a dialectical problem between texts and images. By considering diagrams as hybrid objects, sometimes referred to as iconotexts, a continuum is revealed between these two categories that seemed *a priori* distinctly defined. From then on, diagrams become transmission belts that allow the analysis tools associated with texts to be transported to images and vice versa. In this context, a semiological approach based on the categories introduced by Charles Peirce proves extremely fruitful. As a reminder, Peirce distinguishes iconic signs that have an analogical relationship with what they designate from symbolic signs that do not. Finally, Peirce also introduces indexical signs which physically share the qualities of what they represent:

> A regular progression of one, two, three may be remarked in the three orders of signs, Icon, Index, Symbol. The Icon has no dynamical connection with the object it represents; it simply happens that its qualities resemble those of that object, and excite analogous sensations in the mind for which it is a likeness. But it really stands unconnected to them. The index is physically connected with its object; they make an organic pair, but the interpreting mind has nothing to do with this connection, except remarking it, after it is established. The symbol is connected with its object by virtue of the idea of the symbol-using mind, without which no such connection would exist.[14]

At first sight, images are icons and texts are symbols. But Peirce's semiological theory does not allow for such an identification, since nothing is ever an icon, symbol or index in itself. The nature of objects has no autonomous meaning and must necessarily be complemented by a semiological context consisting of the sign and the interpreter. Together with the object, they therefore form a semiotic triad that must always be considered as a whole, and whose constituents cannot be determined independently of the relationships they have with the rest of the triad[15]. Thus, instead of referring to diagrams as iconotexts, the functional attitude of Peirce's theory suggests that we should call them pictoscriptural signs.
By denying an ontological and autonomous determination to objects, and by inscribing this determination in a relational network that couples an inner psychological state and a sensible phenomenon, Peirce's semiological approach is remarkably well suited to a phenomenological and synthetic approach to diagrammatic practice. But the fact that an object has no autonomous semiotic determination does not prevent it from being isolated mediately through the existence of a bundle of relations that connect different semiotic triads. Peirce considers these relations to be susceptible to formal analysis, from which a true logic of signs can be constructed[16]. The contribution of Peirce's semiology is therefore also consistent with a formalist and analytical attitude to the semantics of diagrams.
These elements, which are not new, have been the subject of abundant commentary, but they shed light on the interest of using Peirce's semiology in the study and editing of diagrams. But the specificity of tables remains a blind spot in this literature, within which the Euclidean diagram is the paradigmatic example of diagrams. Indeed, the table differs from the geometric diagram in several important aspects. Firstly, while the table is an arrangement of objects, which the geometric diagram is also, it presents a symbolic arrangement, where the geometric arrangement is iconic or even indexical. From this point of view, the tabular structure is thus the arbitrary sign of the epistemological structure of the object. However, the young Leibniz's tables show that the geometric, and therefore iconic, aspects of this representation are the witnesses of a signifying practice that should not be masked by historiographic exegesis. Thus, the table is the product of an interaction between two semiological points of view, one purely symbolic, and the other partially iconic. Thus, seeing the young philosopher's tabular practice as a semiotic network allows us to account for both formal structures and instrumental practices.

---

[10] William J. T. Mitchell: "Diagrammatology" in: *Critical Inquiry* 7/3 (1981): pp.622-633.

[11] Frederik Stjernfelt: *Diagrammatology: An Investigation on the Borderlines of Phenomenology, Ontology, and Semiotics.* Springer Science & Business Media, 2007.

[12] Kenneth Manders: "The Euclidean Diagram" in: Paolo Mancosu (ed.): *The Philosophy of Mathematical Practice*, Oxford University Press, 2008, pp.80-133.

[13] Edwin Hutchins: "Material Anchors for Conceptual Blends" in: *Journal of Pragmatics*, Conceptual Blending Theory, 37/10 (2005), pp.1555-1577.

[14] Charles S. Peirce, *Collected Papers Volume 2*, Cambridge, Harvard University Press, 1932, here p.299.

### b. *Michel Serres' structuralism*

The structuralism of Michel Serres constitutes a second extremely fruitful contribution to the study of tables in the young Leibniz. It should be noted, however, that such an appeal requires us to transport his argument to a terrain for which it was not originally intended, and worse, that this deviation of the discourse comes into major conflict with points that are essential to Serres' thinking.
Michel Serres published his doctoral thesis *Le Système de Leibniz et ses modèles mathématiques*[17] in 1968, in which he argues that the idea of structure is the basis for a new approach to Leibnizian thought that takes into account its inherent pluralism without abandoning the project of inscribing it in a single, coherent system. Serres' argument consists in distributing the work of study into two hierarchical levels of analysis and thus making it possible to unify on the upper level points of view that were apparently incompatible on the lower level. This distribution is based on the notion of network structure, of which Michel Serres introduces two substructures, the star and the schema. The star is located at the bottom level and represents the different points of view on the structure. The schema, on the other hand, is located at the upper level and represents the interaction between the different starred points of view on the network. It is thus composed of relations between relations, which leads Serres to assimilate schemas to tables, insofar as they are "series of series". These tables, which Michel Serres describes as "harmonic", have a fundamental property in common: they can be reversed. It is a matter of locally reversing the artificial hierarchy of relations between terms that places intrastellar relations (within stars) below and interstellar relations (between distinct stars) above. This network structure forms the basis of Michel Serres' epistemological structuralism. In an article published two years before his thesis[18], but probably written after it, he explains his idea in a purely abstract way, entirely dissociated from the Leibnizian context in which it was constructed. Serres then expresses the structure as a graph, relying in particular on the image of the chessboard or Penelope's tapestry which gives its name to the article. The network is a "scalene" graph, i.e. any graph that can accommodate all possible sub-graphs. But it is not sub-graphs that are applied to the general structure, but rather dynamic paths, like the stars described in his thesis, whose plasticity produced by the degree of strength of the different internal relations of the graph, presents an evolution over time described by Michel Serres in an organic way.
As we have said, it is not our intention here to apply to our approach this abstract and general network structure, which aims to support a holistic system of thought and to unite worldview and mode of knowledge. In fact, Michel Serres shows in his thesis that this abstract concept derives partly from very concrete elements that he encountered within the Leibnizian corpus on which he relied and in particular the tables. The latter hold a central place in the work of Michel Serres. This is relatively unprecedented in Leibnizian commentary, which had previously focused on the symbolic and formalist practice of the philosopher.
Michel Serres gives a major role to Leibniz's harmonic triangle. His analysis of this subject is mainly based on the account given by the German philosopher in his famous *Historia et Origo Calculi Differentialis*, written in 1714, i.e. more than forty years after most of the facts recounted there. Serres sees this as the paradigmatic example of reversal: by inverting the structure of the table of the arithmetic triangle (i.e. one term of the triangle is the sum of the two terms below it, and not above it), we obtain the table of reversed elements, i.e. inverses. This example thus serves as a quasi-metaphorical archetype for Serres' general thesis. But if we decide to take it seriously and open the file on the genesis of this table (to which Serres did not have access, as these documents were all still unpublished in 1968), we discover a large quantity of Parisian manuscripts filled with numerous tables. The reconstruction of the discovery process reveals a radically different narrative from that on which Serres' analysis is based. In particular, Leibniz did not construct the harmonic triangle by reversing the structure of the arithmetic triangle. But, despite this, the practice of tables is indeed central. And the gradual shift from the "list of lists" strategy to the "double entry table" strategy is a remarkable illustration of the extreme acuity with which, paradoxically, Serres' structural analysis applies. Even on the table present in the *Historia et Origo Calculi Differentialis*, the genetic approach reveals elements that still fit with the properties of Serres' schemas. Indeed, this text exists in two versions, of which only one was published in the 19th century. In both versions of the text, the harmonic triangle does not refer to the same table. In the unpublished version, the triangle consists of the inverses of the terms of the arithmetic triangle, whereas in the edited version, the terms of the triangle are the inverses of the product of the terms of the arithmetic triangle by the rank of these same terms in the triangle.

---

[15] Gérard Deledalle: *Charles S. Peirce's Philosophy of Signs: Essays in Comparative Semiotics*. Indiana University Press, 2000.

[16] Robert Marty: *The Algebra of Signs: An Essay in Scientific Semiotics after Charles Sanders Peirce*, John Benjamins Publishing, 1990.

[17] Michel Serres: *Le Système de Leibniz et ses modèles mathématiques*, Presses universitaires de France, 1968.

[18] Michel Serres: "Penelope ou d'un graphe théorique", in: *Revue Philosophique de la France et de l'Étranger* 156 (1966), pp.41-51.

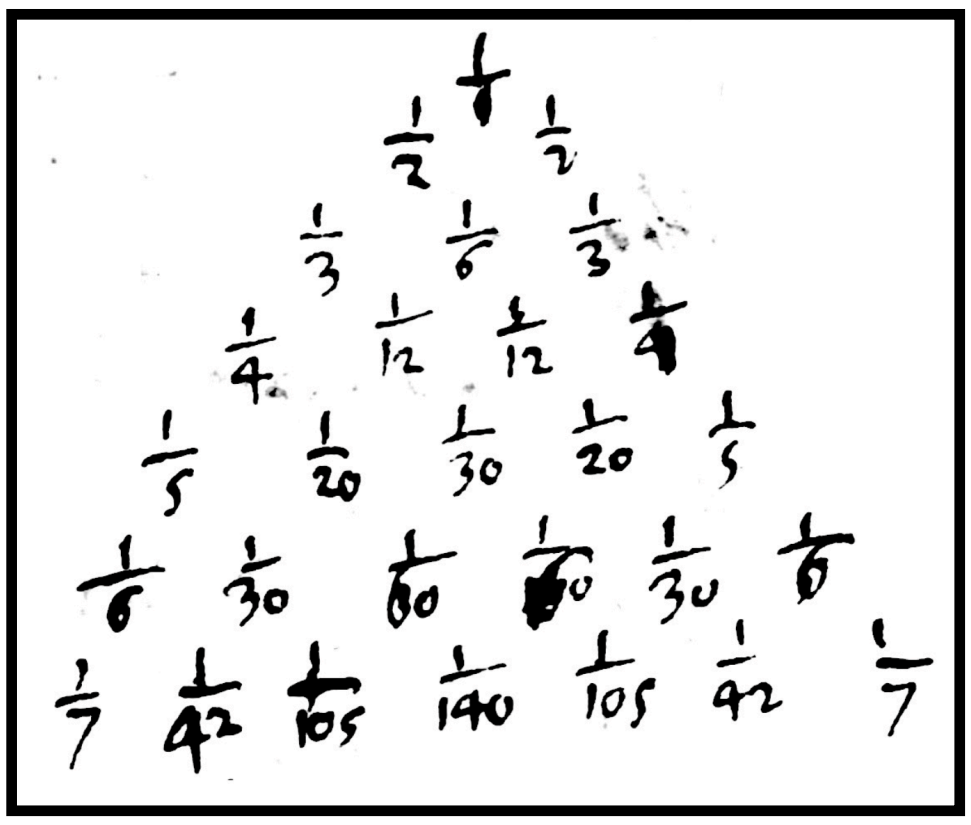
from the edited text

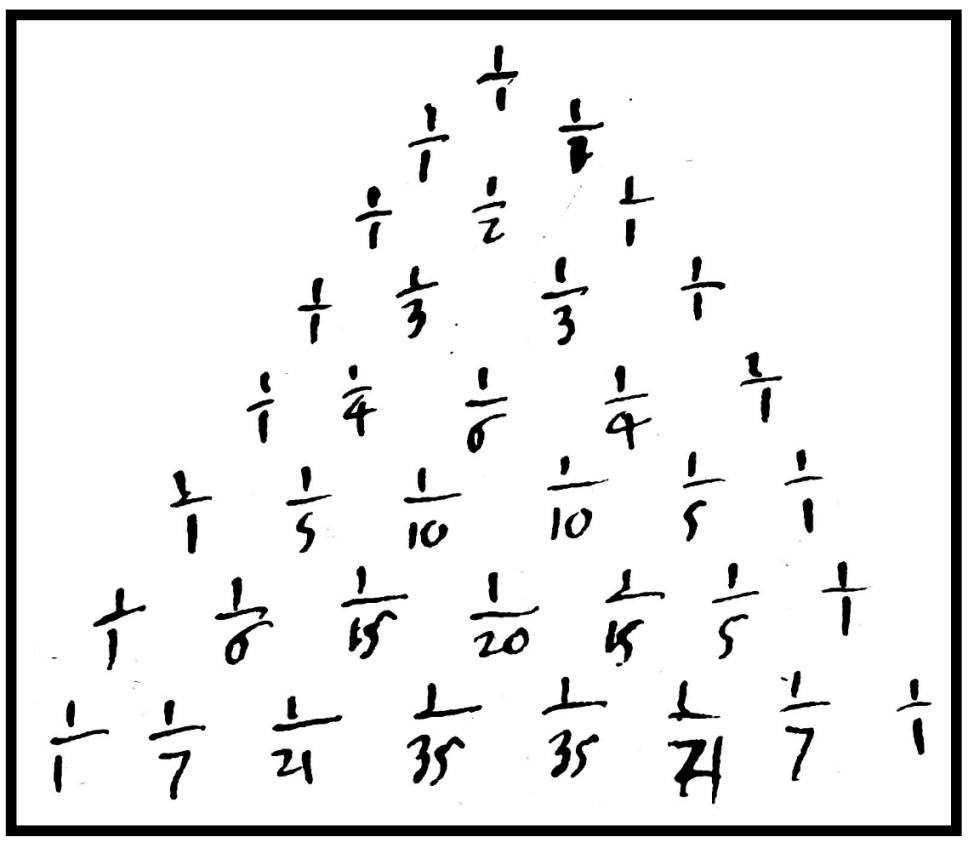
from the unpublished text

In fact, the triangle from the unpublished text appears in the notes of the edited text. Michel Serres therefore had access to both versions, and he chooses to reproduce in his text the version that appears in the notes, i.e. the one from the unpublished text.

| Triangle de Pascal : | | | | | Triangle harmonique : | | | | |
|---|---|---|---|---|---|---|---|---|---|
| 1 | | | | | 1/1 | | | | |
| 1 | 1 | | | | 1/1 | 1/1 | | | |
| 1 | 2 | 1 | | | 1/1 | 1/2 | 1/1 | | |
| 1 | 3 | 3 | 1 | | 1/1 | 1/3 | 1/3 | 1/1 | |
| 1 | 4 | 6 | 4 | 1 | 1/1 | 1/4 | 1/6 | 1/4 | 1/1 |
| etc. | | | | | etc. | | | | |

Inverting the structure of the arithmetic triangle does not give exactly the triangle of inverses, but a slightly modified triangle of inverses, where the terms in each row are divided by the rank of that row. Conversely, the triangle of inverses does not have exactly the inverse structure of the arithmetic triangle, but a slightly modified inverse structure, where instead of adding the two terms below, they must also be multiplied by the ratio between the ranks of the two lines. The harmony that Michel Serres seeks in this example is not perfect, and Leibniz's hesitation between these two versions of the triangle testifies to the fact that he himself was probably well aware of this. But what is most surprising is that Leibniz calls both tables the same: the "harmonic triangle". This means that both tables represent one and the same object. In the first case, one chooses to represent the harmonious structure, but then one has to give up the harmony in the terms, whereas in the second case, it is the terms that are represented harmoniously, but then it is the harmony in the structure that one has to give up. In other words, this clumsy operation, which consists in multiplying the terms of the arithmetic triangle by their rank in the triangle and which breaks the harmony of the inversion, is either implicit in one version (the one chosen by Serres) or explicit in the other version. But, implicit or explicit, this operation is always present within the more general object that Leibniz calls "harmonic triangle". This example thus illustrates once again how two distinct tables can represent the same object. This identification must itself be based on interactions between the two versions, and these interactions can only be described by the genetic approach. Moreover, these interactions are again a wonderful allegorical illustration of Serres' structural theory and how patterns emerge from the entanglement of star paths.
Thus, our approach is born of a reversed movement of Serres', namely to project onto Leibniz' concrete practice elements from his abstract and general conception of structures. This is particularly relevant for the purpose of this contribution, namely to determine a concrete data structure for the digital edition of the young Leibniz's tables that is based on a general philosophical reflection on the nature and function of these artefacts.

## III. Digital implementation: what choices and freedoms?

The main theoretical difficulties in editing tables therefore lie in their profoundly multiform character. Their concrete implementation through the selection of a data structure implies making significant choices: on what criteria can we determine the implicit elements that must be integrated? How should this integration be carried out? However, in addition to these theoretical constraints, there are naturally material and practical constraints. On the one hand, the paper edition being old, it is the repository of past choices in the publishing process. For this reason, the digital edition certainly presents fewer constraints and serves, as this contribution shows, as a support for reflection to future paper editions. On the other hand, digital publishing is not free of limitations, which also need to be taken into account.

### a. *Limitation of computer languages*

The computer language used to produce the paper edition of the complete works of Leibniz is the TeX language, while for the digital edition the TEI language is used. These two languages have a major historical difference: the TeX language was originally developed by Donald Knuth for use in the field of mathematics, whereas the TEI language was developed around scholarly societies specialized in the digital humanities for use in literary studies. Naturally, the TeX language has a large number of tools that are perfectly suited to editing Leibniz's mathematical texts, and even many others have been developed locally in the same language by the editors. The TEI language, on the other hand, has few resources adapted to this purpose. It has therefore been necessary to combine the capabilities of TeX with TEI locally[19], which is facilitated by a property that both languages share but which is also a major problem for critical editing: they are markup languages, i.e. text is always implemented on the basis of a rooted syntax tree whose nodes are delimited by tags. It is then only possible to apply two attributes to two parts of a text if one is included in the other or if they are disjoint. Thus, it is not possible to implement in a formula an attribute which groups together isolated elements which are in separated branches of the syntax tree. This is the case, for example, in the following formula from *Nova de œquationum reductionibus* (A VII, 2, 172):

[19] In fact, we have included the functionalities of the LateX language, derived from TeX, but much more interoperable because it is maintained by an active international community, whereas the TeX language used for print publishing is the product of local maintenance and development.

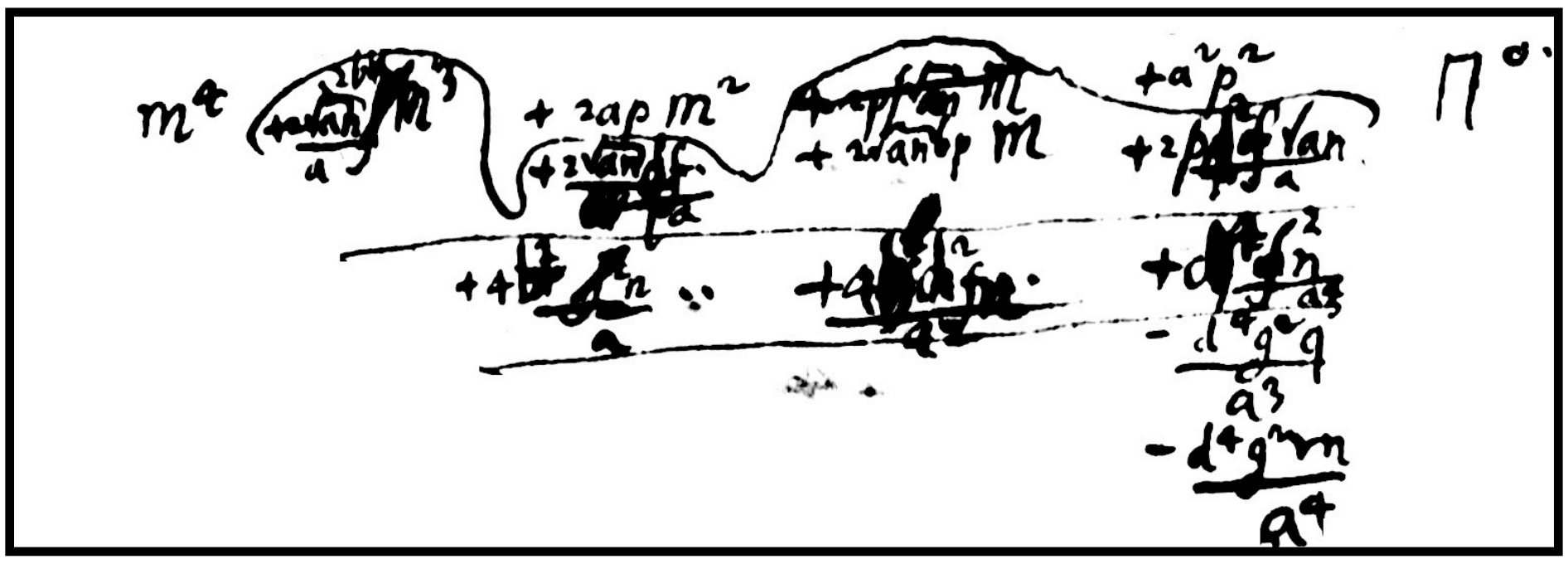

Here, the spatial structure of the term alignments is not compatible with that produced by the drawn lines. The editors, forced to make a choice, chose the structure based on the drawn lines and not on the spatial alignments:

| | | | | | |
|---|---|---|---|---|---|
| $m^4$ | $+$ | $2apm^2$ | | $+$ | $a^2p^2 \sqcap 0.$ |
| $+\dfrac{2\sqrt{an}\,2b}{a}m^3$ | $+$ | $\dfrac{2\sqrt{an}\,df}{a}\,..$ | $+\ 2\sqrt{an}\,bpm$ | $+$ | $\dfrac{2pd^2f\sqrt{an}}{a}$ |
| | $+$ | $\dfrac{4b^2n}{a}\,..$ | $+\ \dfrac{4bd^2fn}{a^2}\,.$ | $+$ | $\dfrac{d^4f^2n}{a^3}$ |
| | | | | $-$ | $\dfrac{d^4g^2q}{a^3}$ |
| | | | | $-$ | $\dfrac{d^4g^2rn}{a^4}$ |

However, this difficulty is constitutive of the tabular structure that we are planning. Indeed, as we have said, the same table must be interpreted in several ways and two distinct tables must possibly be conceived as two different interpretations of the same table.

### *b. Forest structure*

Unlike TeX, for which few solutions have been developed to deal with this problem of syntax conflict, TEI already has tools that have been developed specifically for editorial practices. It is therefore in this language that we must look for a solution. Vincent Buard, a member of the ERC Philiumm project, suggests the forest structure, which is present in the TEI language, as a solution. In graph theory, a forest is only a set of rooted trees. But the forest structure, in the TEI language, makes it possible to enrich this juxtaposition of graphs by applying correspondence attributes that link any two elements of the forest together. Thus, this structure makes it possible to implement formulas with incompatible syntaxes, as in the very simple example below:

[ a + ( b ] - c )

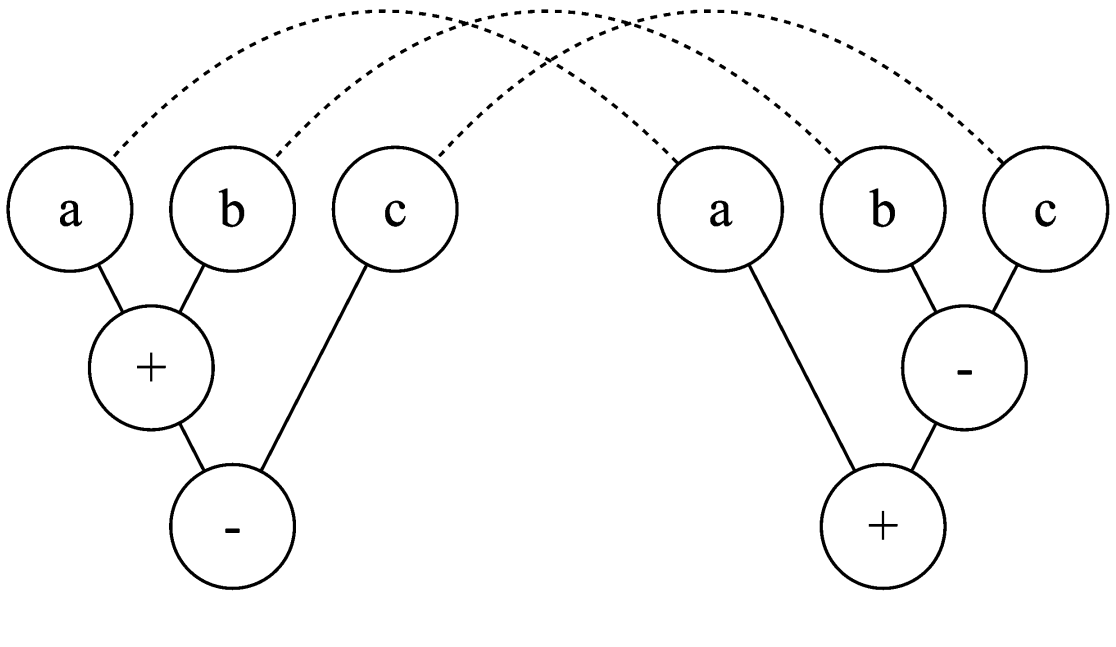

The forest is therefore a data structure, available in the TEI language, which makes it possible to include the reversal property of Michel Serres' schemas. Thus, each tree of the forest corresponds to a possible reading of the table. Moreover, it is possible to add other attributes to the correspondence relations, which are here identification relations, which makes it possible to implement the interactions of the various readings of the table. Moreover, we have shown that these interactions are always based on arguments from a genetic critique of the sources. Thus, the attributes will have to be constructed on the basis of such an approach. Finally, we have shown, thanks to the contribution of Peirce's semiology, that a table always has purely symbolic and partially iconic aspects. Here again, the forest structure of the TEI language allows us to separate these two aspects by including in one tree the symbolic structure and in another the spatial elements which are more of an iconic sign. Let us illustrate this by the very simple example below:

a b

c

In the symbolic structure, we describe the two lines of the table: in the first "a" and "b", and in the second "c". While the iconic structure makes it possible to indicate the fact that "c" is symmetrically related to "a" and "b".

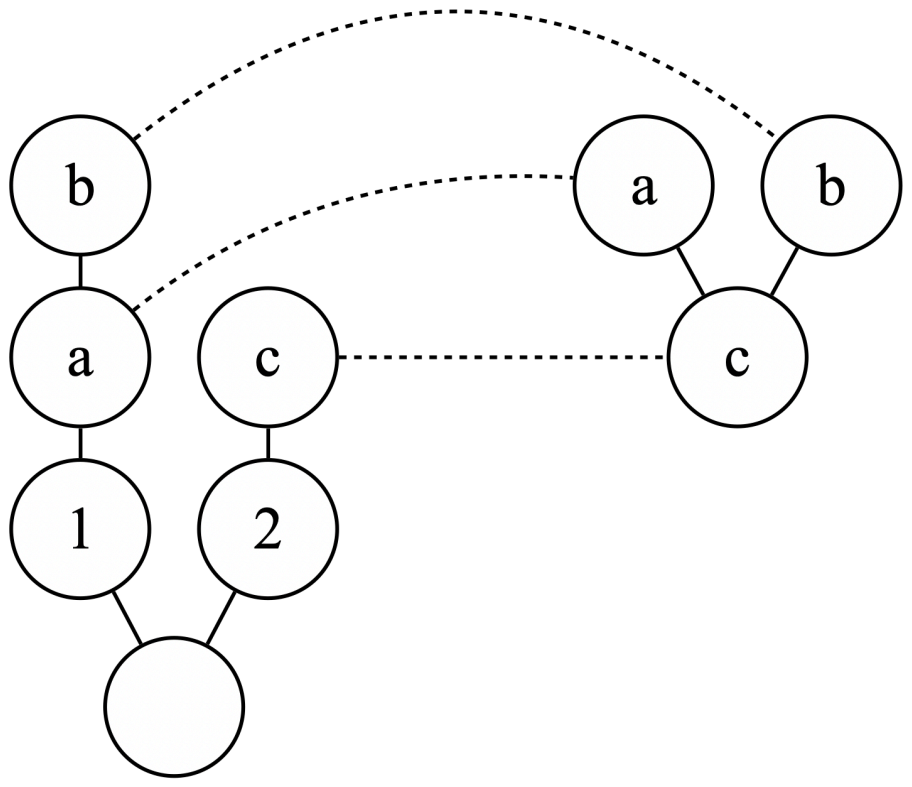

## Conclusion

The tables present in the young Leibniz’s manuscripts are thus very complex philological objects, whose digital edition cannot be satisfied with the classical structure associated with tables. Moreover, a genetic approach to these sources has made it possible to engage in a philosophical and semiological reflection on the nature of these artefacts and raises the question of whether it is relevant or even possible to reduce a table to a linear sequence of signs or instructions. Leibniz's tables are the products of a tabular practice that informs us about the development of his philosophical and practical conception of the art of inventing. But they are also formal writings that are the medium of expression of algorithmic procedures. However, as a symbolic sign, the table always presents an irreducible indeterminacy that allows several procedures to be accommodated on the same support. This concrete indeterminacy in the structure of the particular tables found in the sources fits particularly well with Michel Serres' analysis, which was originally developed to be applied to the abstract and general notion of network structure.

The richness of Leibniz’s tabular structure and practice must be revealed by a genetic approach, which allows, through a kind of stratigraphic analysis, to break with the synchronic conception of diagrams. This contribution of chronology must not serve a psychologism and a naive realism. Michel Serres' notion of a schema, which can accommodate several different starred points of view, makes it possible to account for the multiplicity of possible genetic analyses and their collaborations.

Nevertheless, the technical materiality and the social environment of the tools available to build a digital edition constitute a constraint that pushes us to pragmatism. Thus, without ever losing sight of the impossibility of capturing in a discrete formula the phenomenological complexity of such diagrams, the work of abstraction obliges us to choose the aspects to be kept, simplified or deleted. The forest structure, available in the TEI language, therefore seems to us an excellent candidate: it is an extremely rich structure that allows us to accommodate most of the genetic and semiological aspects that we have presented. Moreover, it is an already existing structure, within a language that is actively maintained by a community dedicated to digital humanity projects.

Moreover, this structure does not only apply to the edition of the tables. In fact, genetic editing involves, in a much more general way, the search for a plastic structure capable of accommodating the various hypotheses resulting from a paleographic analysis. Thus, the forest structure also allows the implementation of equations, as we have already mentioned, but also the implementation of genetic editing of texts for which several versions exist.